\documentclass[suppldata]{interact}
\usepackage{graphicx}
\usepackage{subfigure}
\usepackage{setspace} % to add double space
\usepackage{multirow}
\usepackage{enumerate}
\usepackage{commath, amsmath, amssymb,amsfonts,color,mathtools, enumerate}
\usepackage{hyperref}
\newtheorem{theorem}{Theorem}[section]
\usepackage[linesnumbered,ruled,vlined]{algorithm2e}

\usepackage{moreverb,url}
\usepackage{natbib}% Citation support using natbib.sty
\bibpunct[, ]{(}{)}{;}{a}{}{,}% Citation support using natbib.sty
\renewcommand\bibfont{\fontsize{10}{12}\selectfont}% Bibliography support using natbib.sty

\usepackage{url} % not crucial - just used below for the URL 

\usepackage{longtable}
\usepackage{anyfontsize}
\usepackage{epstopdf}% To incorporate .eps illustrations using PDFLaTeX, etc.

\usepackage{cleveref} % to label eq no to those refered in text
\usepackage{autonum} % to label eq no to those refered in text

\usepackage{afterpage}  % load the afterpage package

\newtheorem*{remark}{Remark}

\begin{document}

\title{Test for symmetry and confidence interval of the parameter $\mu$ of skew-symmetric-Laplace-uniform distribution}

\author{
\name{Raju.~K. Lohot\textsuperscript{a} \href{https://orcid.org/0000-0003-0424-9447}{\includegraphics[scale=0.08]{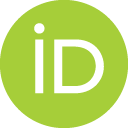}}  \thanks{CONTACT Raju.~K. Lohot. Email: rajulohot.92@gmail.com} and V. U. Dixit\textsuperscript{b} \href{https://orcid.org/0000-0002-4627-589X}{\includegraphics[scale=0.08]{orcid.png}}}
\affil{\textsuperscript{a}Department of Statistics, SVKM's Mithibai College of Arts, Chauhan Institute of Science \& Amrutben Jivanlal College of Commerce and Economics, Vile Parle (W), Mumbai, Maharashtra, India; \textsuperscript{b}Department of Statistics, University of Mumbai, Vidyanagari, Santacruz (E), Mumbai, Maharashtra, India}
}

\maketitle

\begin{abstract}
The skew-symmetric-Laplace-uniform distribution $SSLUD(\mu)$ is introduced in \cite{Lohot2023Dixit} using the skewing mechanism of \cite{azzalini1985class}. Here we derive the most powerful (MP) test for symmetry of the $SSLUD(\mu)$. Since the form of the test statistic is complicated and it is difficult to obtain its exact distribution, critical values and the power of MP test are obtained using simulation. Further, we construct a $100(1-\alpha)\%$ confidence interval (CI) for parameter $\mu$ assuming asymptotic normality and empirical distribution of the maximum likelihood estimator of $\mu$. These two methods are compared based on the average length and coverage probability of the CI. Finally, the CI of the parameter $\mu$ is constructed using data on the ``transformed daily percentage change in the price of NIFTY 50, an Indian stock market index" given in \cite{Lohot2023Dixit}.
\end{abstract}


\begin{amscode}
62F03, 62F25, 65C10
\end{amscode}

\section{Introduction}
Azzalini's form of skew-symmetric density function for any real $\lambda$, as mentioned in \cite{azzalini1985class}, is given as
\begin{equation} \label{azzalini}
    2\,f(x)\,K(\lambda x),
\end{equation} 
where $f$ is a density function symmetric about zero and $K$ is an absolutely continuous distribution function such that the corresponding density function $K\, '$ is symmetric about zero. Using $f$ as the standard Laplace density function and $K$ as a distribution function of Uniform$(-\theta, \theta)$ in (\ref{azzalini}), $SSLUD(\mu)$ is constructed in \cite{Lohot2023Dixit}. The probability density function (pdf) $g(x)$ and cumulative distribution function (cdf) $G(x)$ of $SSLUD(\mu)$ for $\mu \in \mathbb{R}-\{0\}$ are as follows.
\begin{equation}\label{pdf.SSLUD}
    g(x) =
    \begin{cases}
    \; 0 & \text{if} \ \displaystyle \frac{x}{\mu} < -1 ,\\
    \displaystyle\; e^{-\lvert x \rvert} \ \left (\frac{x}{2\mu} + \frac{1}{2} \right)  & \text{if} \ \displaystyle -1 \leqslant \frac{x}{\mu} < 1,  \\
    \displaystyle\; e^{-\lvert x \rvert}\  & \text{if} \ \displaystyle \frac{x}{\mu} \geqslant 1.
    \end{cases} 
\end{equation}

When $\mu<0$,
\begin{subequations}
\begin{equation} \label{SSLUD.neg.cdf}
\begin{split}
G(x) & =
    \begin{cases}
      \displaystyle\; e^x & \text{if} \  x<\mu, \vspace{0.25 cm}\\
      \displaystyle\; \frac{e^{x}}{2\mu}(x+\mu-1)+\frac{e^{\mu}}{2\mu}&  \text{if} \ \mu \leqslant x<0,  \vspace{0.25 cm}\\
     \displaystyle\; 1+\frac{e^{\mu}}{2\mu}-\frac{e^{-x}}{2\mu}(x+\mu+1)&  \text{if} \ 0 \leqslant x< -\mu,  \vspace{0.25 cm}\\
    \; 1 &  \text{if} \  x \geqslant -\mu,
    \end{cases}\\
    \end{split}
    \end{equation}

\noindent \text{and when} $\mu>0$,
\begin{equation}
\label{SSLUD.postive.cdf}
\begin{split}
    G(x) & =
    \begin{cases}
       \; 0 & \text{if} \  x < -\mu, \vspace{0.25 cm} \\
      \displaystyle\; \frac{e^{x}}{2\mu}(x+\mu-1)+\frac{e^{-\mu}}{2\mu}&  \text{if} \ -\mu \leqslant x < 0,  \vspace{0.25 cm}\\
     \displaystyle\; 1+\frac{e^{-\mu}}{2\mu}-\frac{e^{-x}}{2\mu}(x+\mu+1)&  \text{if} \ 0 \leqslant x < \mu,  \vspace{0.25 cm}\\
    \; 1-e^{-x} &  \text{if} \  x \geqslant \mu.
    \end{cases}
\end{split}
\end{equation}
\end{subequations}

\cite{Lohot2023Dixit} provided a comprehensive description of the mathematical properties of (\ref{pdf.SSLUD}) such as moment generating function, characteristic function, first four raw moments, mode, median, hazard rate function, mean deviation about `$a$', R\`enyi entropy, and Shannon entropy. They have discussed the method of maximum likelihood and method of moment for estimating $\mu$. The finite sample and asymptotic properties of these estimators are studied using simulation. Finally, the application of $SSLUD(\mu)$ to real-life data on the ``transformed daily percentage change in the price of NIFTY 50, an Indian stock market index" is discussed. Comparison of fitting of $SSLUD(\mu)$ is done with the fitting of normal distribution $N(\theta, \sigma^2)$, Laplace distribution $L(\theta, \beta)$, and skew-Laplace distribution $SL(\lambda)$ for the above data. 

It can be noted that if $\frac{1}{\mu}=0$ in (\ref{pdf.SSLUD}) then $g(x)$ represents a pdf of standard Laplace distribution which is symmetric about zero. Therefore, we derive the most powerful (MP) test in Section \ref{Sec Test of symmetry} for testing the symmetry of $SSLUD(\mu)$ against positive skewness, that is testing $H_0 : \frac{1}{\mu}=0$ against $H_{1a}: \mu=\mu_1>0$ and similarly, against negative skewness, that is $H_{1b}: \mu=\mu_1<0$. The exact distribution of a test statistic is difficult to obtain; hence, Algorithm \ref{algo:to find k} and Algorithm \ref{algo:to find power} based on simulation are given to obtain critical values and power of the MP test for $H_{1a}$ and $H_{1b}$. Table \ref{tab : cut off point and power} represents critical values and corresponding power of this MP test for different values of $\mu_1\in \mathbb{R}-\{0\}$ and sample size $n=50(50)250$.

Section \ref{Sec. Construction of confidence interval} deals with the confidence interval (CI) for $\mu$. We construct CI for $\mu$ assuming an asymptotic normal distribution of the maximum likelihood estimator (MLE) of $\mu$ using simulation. Also, its modified CI is obtained by removing outliers using the interquartile range (IQR) method of outlier detection. As an alternative method, we construct CI for $\mu$ using the empirical distribution of MLE instead of the normality of MLE. Then modified CI is obtained removing outliers using the same criterion for outlier detection. To compare these two methods of CI construction for $\mu$ assuming asymptotic normality of MLE and using the empirical distribution of MLE, we prepare tables of CI, its average length (AL) and coverage probability (CP) for different values of $\mu$, the sample size $n=50, 100, 150, 250$, and confidence coefficient $1-\alpha=0.95$. Finally in Section \ref{sec application}, using the data on the ``transformed daily percentage change in the price of NIFTY 50, an Indian stock market index" given in \cite{Lohot2023Dixit} which proved to be suitable for $SSLUD(\mu)$, CI for $\mu$  is constructed as an illustration.

\section{Test for symmetry} \label{Sec Test of symmetry}
In this section, we test the hypothesis of symmetry of SSLUD versus positive or negative skewness of the distribution. Since there is no sufficient statistic, the distribution doesn't satisfy the monotone likelihood ratio (MLR) property. First, we obtain the MP test using Neyman Pearson (NP) lemma for testing $H_0 : \frac{1}{\mu}=0$ versus $H_{1a} : \mu=\mu_1, 0<\mu_1<\infty$. The test for alternative $H_{1b} : \mu=\mu_1, -\infty<\mu_1<0$ is obtained on similar lines. As the range of the random variable of $SSLUD(\mu)$ depends on the parameter $\mu$, while using NP lemma we use a joint pdf of order statistics instead of the joint pdf of random sample. The following theorem provides the MP test to test for symmetry against positive skewness in the distribution and its power.

\begin{theorem} \label{thm:MP test for mu>0}
    Suppose $X_1, X_2,\ldots, X_n$ is a random sample which we assume comes from the pdf given in (\ref{pdf.SSLUD}). The MP test of size $\alpha$ for testing $H_0 : \frac{1}{\mu}=0$ versus $H_{1a} : \mu = \mu_1 , 0<\mu_1<\infty$ is
    \begin{equation}
\phi_a (x)  =
\begin{cases}
\; 1   & \text{if} \  \Lambda_a(x)> k_a, \\
\displaystyle\; \gamma_a & \text{if} \ \Lambda_a(x)= k_a, \\
\; 0 & \text{if} \ \Lambda_a(x)< k_a,
\end{cases}
\end{equation}

where \begin{equation}
\Lambda_a(x)  =
\begin{cases}
\; 0   & \text{if} \  y_1 < -\mu_1, \\
\displaystyle\; 2^n \prod_{i=1}^n \left (\frac{y_i}{2\mu_1}+\frac{1}{2}\right)^{I_a(y_i)}   & \text{if} \ -\mu_1 \leqslant y_1 < \mu_1, \\
\; 2^n & \text{if} \ y_1 \geqslant \mu_1,
\end{cases}
\end{equation}
\begin{equation} \label{Indicator a}
    I_a(y_i)=
    \begin{cases}
        \; 1 & \ \text{if} \ -\mu_1 < y_i < \mu_1,\\
        \; 0 & \ otherwise,
    \end{cases}
\end{equation}

\noindent $Y_1, Y_2,\ldots, Y_n$ denote order statistics, and $k_a$ and $\gamma_a$ are such that $E_{H_0}(\phi_a (X))=\alpha$. The power of this MP test is $\beta(\mu_1)=E_{H_{1a}}(\phi_a(X)), \ 0<\mu_1<\infty$. 
\end{theorem}

\begin{proof}

Under $H_0$, the pdf in (\ref{pdf.SSLUD}) is standard Laplace $L(0,1)$ and hence, the joint pdf of the order statistics $Y_1, Y_2,\ldots, Y_n$ is 
    \begin{equation} \label{joint pdf of order stats under H0}
    h_{0}(y_1,\ldots, y_n)  =
        \begin{cases}
            \; \frac{n!}{2^n} e^{-\sum_{i=1}^n |y_i|}   & \text{if} \  -\infty<y_1<y_2<\ldots<y_n<\infty, \\
\displaystyle\; 0  & \ otherwise.
        \end{cases}
    \end{equation}

The joint pdf of the order statistics $Y_1, Y_2,\ldots, Y_n$ under $H_{1a}$ is 
    \begin{equation}
    h_{1a}(y_1,\ldots, y_n)  =
        \begin{cases}
            \; n!\, e^{-\sum_{i=1}^n |y_i|}  \prod_{i=1}^{n} \left (\frac{y_i}{2\mu_1}+\frac{1}{2}\right)^{I_a(y_i)}  & \text{if} \  -\mu_1<y_1<\ldots<y_{r_a}<\\
            & \mu_1<y_{{r_a}+1}<\ldots<y_n<\infty,\\
            & r_a=0, 1,\ldots,n,\\
            \displaystyle\; 0  & \ otherwise,
        \end{cases}
    \end{equation}

\noindent where $r_a$ denotes the number of order statistics lying between $-\mu_1$ and $\mu_1$ provided $y_1 > -\mu_1$, and $I_a(y_i)$ is as defined in (\ref{Indicator a}). Define the ratio $\Lambda_a(x) = \frac{h_{1a}(y_1,\ldots, y_n)}{h_{0}(y_1,\ldots, y_n)}$ which is given as
\begin{equation} \label{lam.a}
\Lambda_a(x)  =
\begin{cases}
\; 0   & \text{if} \  y_1 < -\mu_1, \\
\displaystyle\; 2^n \prod_{i=1}^n \left (\frac{y_i}{2\mu_1}+\frac{1}{2}\right)^{I_a(y_i)}   & \text{if} \ -\mu_1 \leqslant y_1 < \mu_1, \\
\; 2^n & \text{if} \ y_1 \geqslant \mu_1.
\end{cases}
\end{equation}

For $\mu_1 > 0$, the power of MP test $\phi_a(x)$ is 
\begin{equation}
    \beta(\mu_1)=E_{H_{1a}}(\phi_a(X))= P_{H_{1a}}(\Lambda_a(X) > k_a)+\gamma_a \, P_{H_{1a}}(\Lambda_a(X) = k_a).
\end{equation}

Depending on the value of $Y_1$, the computation procedure of cut-off point $k_a$ and $\gamma_a$, and respective power gets divided into three parts as follows. 

\begin{enumerate}
    \item If $y_1 < -\mu_1$ then $\Lambda_a(x) \equiv 0$ that is don't reject $H_0$. Hence, the power $\beta(\mu_1)=E_{H_{1a}}(\phi_a(X))=0$.

    \item If $y_1 > \mu_1$, $P_a = P_{H_0}(Y_1 > \mu_1)= \prod_{i=1}^{n} P_{H_0}(X_i > \mu_1) = \frac{e^{-n\mu_1}}{2^n}$

    \begin{enumerate}
        \item If $\alpha > P_a$ then $\phi_a(x) \equiv 1$ that is reject $H_0$ and its power $\beta(\mu_1)=1$.
        \item If $\alpha \leqslant P_a$ then $\phi_a(x) \equiv \gamma_a$, where $\gamma_a$ is obtained as $\gamma_a = \frac{\alpha}{P_a}$. Its power is $\beta(\mu_1)=\gamma_a$. 
    \end{enumerate}

    \item If $-\mu_1 < y_1 < \mu_1$ then we have to use the following simulation technique to find cut-off point $k_a$. Since the form of $\Lambda_a(x)$ is not analytically tractable, it is very difficult to obtain the exact distribution of $\Lambda_a(x)$. Hence we obtain the cut-off point using Algorithm \ref{algo:to find k} and the power of the test using Algorithm \ref{algo:to find power}.
    
    \begin{algorithm}
        \caption{To find cut-off point $k_a$}
        \label{algo:to find k}
        Generate a random sample of size $n$ from standard Laplace distribution ( i.e. Under $H_0$).\\
        Compute $\Lambda_a(x)$ as given in (\ref{lam.a}) for the sample in Step 1.\\ 
        Repeat Steps 1-2 for $N$ times.\\
        Arrange $N$ values of $\Lambda_a(x)$ in ascending order.\\
        $[N(1-\alpha)]^{th}$ ordered value of $\Lambda_a(x)$ is a value of cut off point $k_a$. Hence, the test is 
            $$
            \phi_a(x)  =
            \begin{cases}
             1   & \text{if} \  \Lambda_{a}(x)> k_a, \\
             0 & \text{if} \ \Lambda_{a}(x) < k_a.
            \end{cases}
            $$
    \end{algorithm}

    \begin{algorithm}
        \caption{To find the power of $\phi_a(x)$ at $\mu_1$, $0 < \mu_1 < \infty$}
        \label{algo:to find power}
        Generate a random sample of size $n$ from $SSLUD(\mu_1)$ ( i.e. Under $H_{1a}$).\\
        Compute $\Lambda_a(x)$ as given in (\ref{lam.a}) for the sample in Step 1.\\ 
        Repeat Steps 1-2 for $N$ times.\\
        Obtain the proportion of $\Lambda_a(x)$ values greater than $k_a$ out of $N$ values obtained in Step 3, where $k_a$ is as obtained in Step 5 of Algorithm \ref{algo:to find k}. This proportion is the power $\beta(\mu_1)$ obtained by using the simulation technique.
        \end{algorithm}
\end{enumerate}

\end{proof}

Following Theorem \ref{thm:MP test for mu<0} provides the MP test to test for symmetry against negative skewness of the distribution and its power.

\begin{theorem} \label{thm:MP test for mu<0}
    Suppose $X_1, X_2,\ldots, X_n$ is a random sample which we assume comes from the pdf given in (\ref{pdf.SSLUD}). The MP test of size $\alpha$ for testing $H_0 : \frac{1}{\mu}=0$ versus $H_{1b} : \mu=\mu_1, -\infty<\mu_1<0$ is

    \begin{equation}
\phi_b (x)  =
\begin{cases}
\; 1   & \text{if} \  \Lambda_b(x)> k_b, \\
\displaystyle\; \gamma_b & \text{if} \ \Lambda_b(x)= k_b, \\
\; 0 & \text{if} \ \Lambda_b(x)< k_b,
\end{cases}
\end{equation}

where \begin{equation}
\Lambda_b(x)  =
\begin{cases}
\; 2^n & \text{if} \ y_n < \mu_1,\\
\displaystyle\; 2^n \prod_{i=1}^n \left (\frac{y_i}{2\mu_1}+\frac{1}{2}\right)^{I_b(y_i)}   & \text{if} \ \mu_1 \leqslant y_n < -\mu_1, \\
\; 0   & \text{if} \  y_n \geqslant -\mu_1,

\end{cases}
\end{equation}
\begin{equation}    \label{indicator b}
    I_b(y_i)=
    \begin{cases}
        \; 1 & \ \text{if} \ \mu_1 < y_i < -\mu_1,\\
        \; 0 & \ otherwise,
    \end{cases}
\end{equation}

\noindent $Y_1, Y_2,\ldots, Y_n$ denote order statistics, and $k_b$ and $\gamma_b$ are such that $E_{H_0}(\phi_b (X))=\alpha$. The power of this MP test is given by $\beta(\mu_1)=E_{H_{1b}}(\phi_b(X)), \ -\infty <\mu_1<0$. 
\end{theorem}
\begin{proof}

The joint pdf of the order statistics $Y_1, Y_2,\ldots, Y_n$ under $H_0$ is as given in (\ref{joint pdf of order stats under H0}) and under $H_{1b}$ is 

    \begin{equation}
    h_{1b}(y_1,\ldots, y_n)  =
        \begin{cases}
            \; n!\, e^{-\sum_{i=1}^n |y_i|}  \prod_{i=1}^{n} \left (\frac{y_i}{2\mu_1}+\frac{1}{2}\right)^{I_b(y_i)}  & \text{if} \  -\infty<y_1<\ldots<y_{n - r_b}<\\
            & \mu_1<y_{n - r_b+1}<\ldots<y_n<-\mu_1,\\ & r_b=0, 1,\ldots,n,\\
            \displaystyle\; 0  & \ otherwise,
        \end{cases}
    \end{equation}

\noindent where $r_b$ denotes the number of order statistics lying between $\mu_1$ and $-\mu_1$ provided $y_n < -\mu_1$, and $ I_b(y_i)$ is as defined in (\ref{indicator b}). Define the ratio $\Lambda_b(x) = \frac{h_{1b}(y_1,\ldots, y_n)}{h_{0}(y_1,\ldots, y_n)}$ which is given as
\begin{equation} \label{lam.b}
\Lambda_b(x)  =
\begin{cases}
\; 2^n & \text{if} \ y_n < \mu_1,\\
\displaystyle\; 2^n \prod_{i=1}^n \left (\frac{y_i}{2\mu_1}+\frac{1}{2}\right)^{I_b(y_i)}   & \text{if} \ \mu_1 \leqslant y_n < -\mu_1, \\
\; 0   & \text{if} \  y_n \geqslant -\mu_1.
\end{cases}
\end{equation}

For $\mu_1 < 0$, the power of MP test $\phi_b(x)$ is 
$$
\beta(\mu_1)=E_{H_{1b}}(\phi_b(X))= P_{H_{1b}}(\Lambda_b(X) > k_b)+\gamma_b \, P_{H_{1b}}(\Lambda_b(X) = k_b).    
$$

Depending on the value of $Y_n$, the computation procedure of cut-off point $k_b$ and $\gamma_b$, and respective power gets divided into three parts as follows.

\begin{enumerate}
    \item If $y_n > -\mu_1$, then $\Lambda_b(x) \equiv 0$ that is don't reject $H_0$. Hence, the power $\beta(\mu_1)=0$.

    \item If $y_n < \mu_1$, $P_b = P_{H_0}(Y_n < \mu_1)= \prod_{i=1}^{n} P_{H_0}(X_i < \mu_1) = \left[ e^{\mu_1} \, (1-\frac{e^{\mu_1}}{2}) \right ]^n $.

    \begin{enumerate}
        \item If $\alpha > P_b$, then $\phi_b(x) \equiv 1$ that is reject $H_0$ and its power $\beta(\mu_1)=1$.
        \item If $\alpha \leqslant P_b$, then $\phi_b(x) \equiv \gamma_b$, where $\gamma_b$ is obtained as $\gamma_b = \frac{\alpha}{P_b}$. Its power is $\beta(\mu_1)=\gamma_b$. 
    \end{enumerate}

    \item If $\mu_1 < y_n < -\mu_1$, then cut-off point is obtained using the simulation procedure given in Algorithm \ref{algo:to find k} by replacing $\Lambda_a(x)$ with $\Lambda_b(x)$. Here the MP test $\phi_b(x)$ is 

     \begin{equation}
            \phi_b(x)  =
            \begin{cases}
            \; 1   & \text{if} \  \Lambda_{b}(x)> k_b, \\
            \; 0 & \text{if} \ \Lambda_{b}(x)\leqslant k_b.
            \end{cases}
    \end{equation}
    The power $\beta(\mu_1)$ of this MP test of size $\alpha$ can be obtained using Algorithm \ref{algo:to find power} by replacing $\Lambda_a(x)$, $H_{1a}$, and $k_a$ with $\Lambda_b(x)$, $H_{1b}$, and $k_b$ respectively.
    
\end{enumerate}
\end{proof}

As an illustration, the values of cut-off point $k$ and power of this MP test of size $\alpha$ are calculated for sample sizes 50(50)250 for testing $H_0: \frac{1}{\mu}=0$ against $H_{1a}: \mu = \mu_1, 0 < \mu_1 < \infty$ for $\mu_1 = 1, 2, 3, 4, 5, 10, 15, 25$. Similarly, cut off point and power for testing $H_0: \frac{1}{\mu}=0$ against $H_{1b}: \mu = \mu_1, -\infty < \mu_1 < 0$ for $\mu_1 = -25, -15, -10, -5, -4, -3, -2, -1$ are calculated. These values are presented in Table \ref{tab : cut off point and power}. We have taken the simulation size as $N=1000$. Numbers in parenthesis represent the value of cut-off point $k$ and without parenthesis shows the power of the test.

\begin{table} 
\caption{The values of cut-off point $k$ and power of the MP test of size $\alpha=0.05$ for sample size n for testing $H_0: \frac{1}{\mu}=0$ against $H_{1}: \mu = \mu_1$ using simulation size $N=1000$} 
\label{tab : cut off point and power}
\begin{center}
\renewcommand{\arraystretch}{1.2}
\begin{tabular}{cccccc} 
\hline \multirow{2}{1 cm}{$\mu_1$} & \multicolumn{5}{c}{$n$}\\ \cline{2-6}
& 50 & 100 & 150 & 200 & 250 \\ \hline 

\multirow{2}{1 cm}{-25} & 0.095	& 0.136	& 0.175	& 0.184	& 0.236\\
 & (1.81046)	& (2.13942)	& (2.43615)	& (2.86453)	& (2.78399)\\
\hline 
\multirow{2}{1 cm}{-15} & 0.159	& 0.26	& 0.304	& 0.365	& 0.467\\
& (2.36273)	& (2.77728)	& (3.83689)	& (3.80542)	& (3.53951)\\
\hline 
\multirow{2}{1 cm}{-10} & 0.29	& 0.422	& 0.556	& 0.688	& 0.709\\
& (3.03149)	& (4.16676)	& (3.94934)	& (3.12039)	& (3.4524)\\
\hline 
\multirow{2}{1 cm}{-5} & 0.662	& 0.923	& 0.984	& 0.996	& 1\\
& (4.26374)	& (1.34236)	& (0.63367)	& (0.0763)	& (0.0129)\\
\hline 
\multirow{2}{1 cm}{-4} & 0.868	& 0.985	& 0.999	& 1	& 1\\
& (2.46656)	& (0.53052)	& (0.01091)	& (0.00016)	& (0)\\
\hline 
\multirow{2}{1 cm}{-3} & 0.984	& 1	& 1	& 1 & 1\\
& (0.6219)	& (0.00013)	& (0)	& (0)	& (0)\\
\hline 
\multirow{2}{1 cm}{-2} & 1	& 1	& 1	& 1	& 1\\
& (0)	& (0)	& (0)	& (0)	& (0)\\
\hline 
\multirow{2}{1 cm}{-1} & 1	& 1	& 1	& 1	& 1\\
& (0)	& (0)	& (0)	& (0)	& (0)\\
\hline 
\multirow{2}{1 cm}{1} & 1	& 1	& 1	& 1	& 1\\
& (0)	& (0)	& (0)	& (0)	& (0)\\
\hline 
\multirow{2}{1 cm}{2} & 1	& 1	& 1	& 1	& 1\\
& (0)	& (0)	& (0)	& (0)	& (0)\\
\hline 
\multirow{2}{1 cm}{3} & 0.982	& 1	& 1	& 1	& 1\\
& (0.79937)	& (0.00001)	& (0)	& (0)	& (0)\\
\hline 
\multirow{2}{1 cm}{4} & 0.855	& 0.987	& 1	& 1	& 1\\
& (2.49584)	& (0.3491)	& (0.01026)	& (0.00006)	& (0) \\
\hline 
\multirow{2}{1 cm}{5} & 0.673	& 0.928	& 0.981	& 0.997	& 0.999\\
& (3.65291)	& (1.54359)	& (0.61227)	& (0.0741)	& (0.0149)\\
\hline 
\multirow{2}{1 cm}{10} & 0.321	& 0.414	& 0.524	& 0.667	& 0.732\\
& (2.76303)	& (3.85758)	& (4.71173)	& (3.05039)	& (3.37577)\\
\hline 
\multirow{2}{1 cm}{15} & 0.165	& 0.245	& 0.401	& 0.396	& 0.433\\
& (2.39985)	& (3.29963)	& (2.82202)	& (3.74806)	& (3.83881)\\
\hline 
\multirow{2}{1 cm}{25} & 0.117	& 0.153	& 0.188	& 0.198	& 0.204\\
& (1.77401)	& (2.11612)	& (2.40579)	& (2.72554)	& (2.96513)\\
\hline 
\end{tabular}
\label{Table of Bias and MSE 2}
\end{center}
\end{table}

Since the values of cut-off point $k$ depend on the actual value of $\mu_1$, the same test can not be the UMP test for the alternative hypothesis $H_{1a}: 0<\mu<\infty$ or $H_{1b}: -\infty<\mu<0$. Therefore, the UMP test does not exist. From Table \ref{tab : cut off point and power}, it is seen that for fixed $\mu_1$, the power increases as the sample size $n$ increases as expected because the MP test is consistent. Similarly, for a fixed value of $n$, the power increases as $\mu_1$ approaches to zero.
\begin{remark}
    In part 2(b) of proof of theorem \ref{thm:MP test for mu>0}, to find the values of sample size $n$ for which $\alpha < P_a$ so that $\phi_a(x) \equiv \gamma_a$, we have the inequality as follows. 
    \begin{equation}
        n < \frac{-\, ln\alpha}{\mu_1+ln2}
    \end{equation}
    Similarly, in part 2(b) of proof of theorem \ref{thm:MP test for mu<0}, to find the values of sample size $n$ for which $\alpha < P_b$ so that $\phi_b(x) \equiv \gamma_b$, we have the inequality as follows. 
    \begin{equation}
        n < \frac{ln\alpha}{\mu_1 + ln\left(1-\frac{e^{\mu_1}}{2}\right)}
    \end{equation}

\begin{table} 
\caption{Maximum value of $n$ for $\alpha = 0.01, 0.05$ and $\mu_1$ such that $\alpha < P_a$ for $\mu_1 > 0$ and $\alpha < P_b$ for $\mu_1 < 0$}
\label{tab: max n}
\begin{center}
\renewcommand{\arraystretch}{1.5}
{\begin{tabular}{lcccccc} 
\hline \multirow{2}{1 cm}{$\alpha$} & \multicolumn{6}{c}{$\mu_1$}\\ \cline{2-7}

& -5 & -3 & -1 & 1 & 3 & 5 \\
\hline
$0.01$ & 6 & 4 & 3 & 3 & 2 & 2\\
 \hline
0.05 & 4 & 3 & 2 & 2 & 2 & 1\\
\hline
\end{tabular}}
\label{Table of AIC BIC 2}
\end{center}
\end{table}

\noindent Table \ref{tab: max n} represents maximum value of the sample size $n$ for $\alpha=0.01, 0.05$ and $\mu_1 =-5, -3, -1, 1, 3, 5$ such that $\alpha < P_a$ for $\mu_1 > 0$ and $\alpha < P_b$ for $\mu_1 < 0$.  
\end{remark}

\section{Construction of confidence interval for \texorpdfstring{$\mu$}{TEXT}}
\label{Sec. Construction of confidence interval}

Given an iid sample $X_1, X_2,\ldots,X_n$ which we assume coming from the $SSLUD(\mu)$, CI can be constructed for the unknown parameter $\mu$ by using one of the following two methods based on simulation. To assess the accuracy of CI, here we compute AL and CP of the CI for different values of $\mu$ and sample size $n$. Each of these two methods is further modified to increase accuracy by using the IQR method of outlier detection.

\subsection{Assuming an asymptotic normal distribution of the MLE \texorpdfstring{$\hat{\mu}$}{TEXT}}
\label{Sec. construction of CI using asymptotic normal distribution of the MLE}
It is obvious that $SSLUD(\mu)$ given in (\ref{pdf.SSLUD}) does not satisfy the regularity conditions and the MLE $\hat{\mu}$ is not in closed form but, obtained by maximizing the likelihood function numerically. Hence here we construct $100(1-\alpha)\%$ CI for $\mu$ assuming asymptotic normality of the MLE $\hat{\mu}$ as given in Algorithm \ref{algo:CI using asymp normality of MLE}. To validate the procedure of construction of CI for $\mu$ given in Algorithm \ref{algo:CI using asymp normality of MLE}, we compute the AL and CP of the CI using simulation which is explained in Algorithm \ref{algo:AL and CP for CI using asymp normality of MLE}. 

\begin{algorithm}
    \caption{To compute CI for $\mu$ assuming asymptotic normality of MLE $\hat{\mu}$}
    \label{algo:CI using asymp normality of MLE}
    Compute the MLE $\hat{\mu}$ of $\mu$ for given iid sample $x_1, x_2, \ldots, x_n$.\\
    Generate a random sample of size $n$ from $SSLUD(\hat{\mu})$ and obtain its MLE, say $\hat{\mu}_1$.\\
    Repeat Step 2, $N$ times. The MLEs obtained using these $N$ samples are denoted by $\hat{\mu}_i$, $i=1, 2, \ldots, N$.\\
    Obtain mean and variance of these $N$ values of MLE. Consider this variance as the estimate of the variance of the MLE $\hat{\mu}$, say $\hat{V}(\hat{\mu})$.\\
    $100(1-\alpha)\%$ CI for $\mu$ is $\hat{\mu}\pm z_{\alpha/2} \sqrt{\hat{V}(\hat{\mu})}$, where $z_{\alpha/2}$ denotes upper $({\alpha/2})^{th}$ quantile of standard normal distribution.
\end{algorithm}

\begin{algorithm}
    \caption{To compute the AL and CP of the CI for $\mu$ constructed assuming asymptotic normality of MLE $\hat{\mu}$}
    \label{algo:AL and CP for CI using asymp normality of MLE}
    Generate a random sample of size $n$ from $SSLUD(\mu)$ for a specified value of $\mu$.\\
    For the random sample generated in Step 1 follow Steps 1-4 of Algorithm \ref{algo:CI using asymp normality of MLE} to compute $\hat{V}(\hat{\mu})$.\\
    Go to Step 1. Compute MLE of $\mu$ say $\hat{\mu_1}$ and CI as $\hat{\mu}_1\pm z_{\alpha/2} \sqrt{\hat{V}(\hat{\mu})}$.\\
    Repeat Step 3, $N$ times such that $N$ CIs for $\mu$ are available as $\hat{\mu}_i\pm z_{\alpha/2} \sqrt{\hat{V}(\hat{\mu})}$, $i=1, 2, \ldots, N$.\\
    AL of these CIs is $2 z_{\alpha/2} \sqrt{\hat{V}(\hat{\mu})}$ and its CP is the proportion of CIs obtained in Step 4 which contain $\mu$ specified in Step 1.
\end{algorithm}

In the following we construct the modified CI by removing the outliers from the $N$ values of $\hat{\mu}_i$, $i=1, 2, \ldots, N$ obtained in Step 3 of Algorithm \ref{algo:CI using asymp normality of MLE} to make this CI robust concerning outliers. Suppose $Q_1$ and $Q_3$ denote the $1^{st}$ and $3^{rd}$ quartile of $N$ values of $\hat{\mu}_i$ respectively. Here, to remove outliers, we use the IQR method of outlier detection. Any value $\hat{\mu}_i$, $i=1, 2, \ldots, N$ which is less than $Q_1 -1.5 IQR$ or greater than $Q_3+1.5 IQR$ is removed and based on remaining $N^*$ values of $\hat{\mu}_i$, $\hat{V}(\hat{\mu})$ is computed in Step 4 of Algorithm \ref{algo:CI using asymp normality of MLE}. Thus, the modified CI is constructed using Step 5 of Algorithm \ref{algo:CI using asymp normality of MLE}, and its AL and CP is computed using Algorithm \ref{algo:CI using asymp normality of MLE} with $\hat{V}(\hat{\mu})$ computed as discussed above.

Table \ref{Table of CI Normality 1}, Table \ref{Table of CI Normality 2}, and Table \ref{Table of CI Normality 3} represents $\hat{V}(\hat{\mu})$, CI of $\mu$ assuming asymptotic normality of MLE $\hat{\mu}$, its AL and CP, and all these entities removing outliers for $\mu=-5, -4, -3, -2, -1.5, -1, -0.75, -0.5, -0.25, 0.25, 0.5, 0.75, 1, 1.5, 2, 3, 4, 5$ and the sample size $n= 50, 100, 150, 250$ using simulation size $N=1000$. In each cell, the first value presents $\hat{V}(\hat{\mu})$, values in the parenthesis represent CI for $\mu$, and the next two values present AL and CP of CI respectively. Similarly, in each cell values in bold type present modified $\hat{V}(\hat{\mu})$, modified CI for $\mu$, and AL and CP of modified CI respectively.

From these tables, it is seen that the CI constructed assuming the normality of the MLE $\hat{\mu}$ works well concerning the CP. Further when $\hat{V}(\hat{\mu})$ is large due to outliers, a modified method of construction of CI performs better concerning $\hat{V}(\hat{\mu})$ and AL. Also, it can be noted that as the sample size $n$ increases or the value of $\mu$ approaches zero, $\hat{V}(\hat{\mu})$ decreases with few exceptions.

\begin{table}
\caption{Estimate of the variance of MLE \texorpdfstring{$\hat{\mu}$}{TEXT}, CI of $\mu$ assuming an asymptotic normality of MLE \texorpdfstring{$\hat{\mu}$}{TEXT}, its AL and CP, and all these entities removing outliers for $\mu$, sample size $n$, and simulation size $N=1000$}
\renewcommand{\arraystretch}{1}
\begin{tabular*}{\textwidth}{@{\extracolsep\fill}ccccc} 
\toprule
\multirow{2}{1 cm}{$\mu$} & \multicolumn{4}{@{}c@{}}{n}  \\ \cmidrule{2-5}%
& 50 & 100 & 150 & 250 \\ \midrule
\multirow{6}{1 cm}{-5} & 9030.238  & 249.492 & 1.573 & 4.196\\
 & (-193.195, 179.307) & (-36.167, 25.749) & (-7.064, -2.148) & (-9.963, -1.933)\\
 & 372.501, 0.998 & 61.916, 0.996 & 4.916, 0.931 & 8.029, 0.996 \\
 & \textbf{10.297} & \textbf{2.119} & \textbf{0.855} & \textbf{1.456 }\\
 & \textbf{(-13.233, -0.655)} & \textbf{(-8.062, -2.356)} & \textbf{(-6.419, -2.793)} & \textbf{(-8.313, -3.583)} \\
 & \textbf{12.579, 0.9} & \textbf{5.706, 0.887} & \textbf{3.625, 0.874} & \textbf{4.731, 0.968} \\
\bottomrule
\multirow{6}{1 cm}{-4} &  0.144 & 0.464 & 113.102 & 0.357\\
 & (-2.526, -1.038) & (-4.538, -1.866) & (-26.79, 14.898) & (-5.062, -2.718)\\
 & 1.489, 0.406 & 2.671, 0.838 & 41.688, 1 & 2.344, 0.953 \\
 & \textbf{0.135} & \textbf{0.349} & \textbf{2.666} & \textbf{0.271}\\
 & \textbf{(-2.503, -1.061)} & \textbf{(-4.361, -2.043)} & \textbf{(-9.146, -2.746)} & \textbf{(-4.91, -2.87)} \\
 & \textbf{1.441, 0.387} & \textbf{2.317, 0.774} & \textbf{6.4, 0.993} & \textbf{2.04, 0.919} \\
\bottomrule
\multirow{6}{1 cm}{-3} &  4.72 & 0.301 & 0.234 & 0.114\\
 & (-7.438, 1.078) & (-3.794, -1.644) & (-4.006, -2.11) & (-3.625, -2.299)\\
 & 8.516, 0.991 & 2.151, 0.942 & 1.896, 0.96 & 1.326, 0.943 \\
 & \textbf{0.674} & \textbf{0.193} & \textbf{0.198} & \textbf{0.103}\\
 & \textbf{(-4.789, -1.571)} & \textbf{(-3.581, -1.857)} & \textbf{(-3.93, -2.186)} & \textbf{(-3.591, -2.333)} \\
 & \textbf{3.218, 0.931} & \textbf{1.723, 0.892} & \textbf{1.744, 0.945} & \textbf{1.259, 0.929} \\
\bottomrule
\multirow{6}{1 cm}{-2} & 0.174  & 0.046 & 0.057 & 0.021\\
 & (-2.711, -1.077) & (-1.895, -1.059) & (-2.451, -1.515) & (-2.004, -1.434)\\
 & 1.635, 0.925 & 0.836, 0.82 & 0.936, 0.933 & 0.571, 0.876 \\
 & \textbf{0.143} & \textbf{0.042} & \textbf{0.051} & \textbf{0.019}\\
 & \textbf{(-2.636, -1.152)} & \textbf{(-1.88, -1.074)} & \textbf{(-2.426, -1.54)} & \textbf{(-1.99, -1.448)} \\
 & \textbf{1.484, 0.889} & \textbf{0.806, 0.808} & \textbf{0.885, 0.922} & \textbf{0.541, 0.858} \\
\bottomrule
\multirow{6}{1 cm}{-1.5} & 0.099 & 0.04 & 0.026 & 0.013 \\
 & (-2.133, -0.899) & (-1.851, -1.069) & (-1.74, -1.112) & (-1.605, -1.167)\\
 & 1.234, 0.935 & 0.782, 0.932 & 0.629, 0.946 & 0.438, 0.909\\
 & \textbf{0.093} & \textbf{0.036} & \textbf{0.023} & \textbf{0.011}\\
 & \textbf{(-2.114, -0.918)} & \textbf{(-1.832, -1.088)} & \textbf{(-1.724, -1.128)} & \textbf{(-1.592, -1.18)} \\
 & \textbf{1.196, 0.927} & \textbf{0.744, 0.909} & \textbf{0.597, 0.931} & \textbf{0.412, 0.884} \\
\bottomrule
\multirow{6}{1 cm}{-1} &  0.062 & 0.018 & 0.01 & 0.006\\
 & (-1.654, -0.676) & (-1.227, -0.699) & (-1.024, -0.634) & (-1.053, -0.753)\\
 & 0.977, 0.974 & 0.528, 0.927 & 0.39, 0.906 & 0.3, 0.922 \\
 & \textbf{0.055} & \textbf{0.016} & \textbf{0.009} & \textbf{0.006}\\
 & \textbf{(-1.624, -0.706)} & \textbf{(-1.209, -0.717)} & \textbf{(-1.018, -0.64)} & \textbf{(-1.049, -0.757)} \\
 & \textbf{0.918, 0.963} & \textbf{0.492, 0.907} & \textbf{0.378, 0.898} & \textbf{0.292, 0.915} \\
\bottomrule
\end{tabular*}
\label{Table of CI Normality 1}
\end{table}

\begin{table}
\caption{Estimate of the variance of MLE \texorpdfstring{$\hat{\mu}$}{TEXT}, CI of $\mu$ assuming an asymptotic normality of MLE \texorpdfstring{$\hat{\mu}$}{TEXT}, its AL and CP, and all these entities removing outliers for $\mu$, sample size $n$, and simulation size $N=1000$}
\renewcommand{\arraystretch}{1}
\begin{tabular*}{\textwidth}{@{\extracolsep\fill}ccccc} 
\toprule
\multirow{2}{1 cm}{$\mu$} & \multicolumn{4}{@{}c@{}}{n}  \\ \cmidrule{2-5}%
& 50 & 100 & 150 & 250 \\ \midrule
\multirow{6}{1 cm}{-0.75} & 0.029 & 0.01 & 0.008 & 0.004\\
 & (-1.115, -0.451) & (-0.844, -0.45) & (-0.913, -0.565) & (-0.834, -0.59)\\
 & 0.665, 0.936 & 0.395, 0.907 & 0.349, 0.937 & 0.245, 0.926 \\
 & \textbf{0.026} & \textbf{0.01} & \textbf{0.007} & \textbf{0.004}\\
 & \textbf{(-1.099, -0.467)} & \textbf{(-0.839, -0.455)} & \textbf{(-0.907, -0.571)} & \textbf{(-0.829, -0.595)} \\
 & \textbf{0.631, 0.923} & \textbf{0.384, 0.901} & \textbf{0.337, 0.931} & \textbf{0.233, 0.915} \\
\bottomrule
\multirow{6}{1 cm}{-0.5} &  0.009 & 0.008 & 0.003 & 0.002\\
 & (-0.486, -0.11) & (-0.723, -0.373) & (-0.529, -0.303) & (-0.528, -0.352)\\
 & 0.377, 0.832 & 0.35, 0.957 & 0.227, 0.886 & 0.177, 0.897 \\
 & \textbf{0.009} & \textbf{0.007} & \textbf{0.003} & \textbf{0.002}\\
 & \textbf{(-0.486, -0.11)} & \textbf{(-0.711, -0.385)} & \textbf{(-0.524, -0.308)} & \textbf{(-0.523, -0.357)} \\
 & \textbf{0.375, 0.832} & \textbf{0.326, 0.941} & \textbf{0.217, 0.874} & \textbf{0.165, 0.886} \\
\bottomrule
\multirow{6}{1 cm}{-0.25} & 0.001 & 0.004 &  0.002 & 0.001 \\
 & (-0.143, 0.003) & (-0.378, -0.146) & (-0.332, -0.168) & (-0.34, -0.21)\\
 & 0.146, 0.597 & 0.233, 0.918 & 0.163, 0.906 & 0.13, 0.929 \\
 & \textbf{0.001} & \textbf{0.003} & \textbf{0.002} & \textbf{0.001}\\
 & \textbf{(-0.14, 0)} & \textbf{(-0.368, -0.156)} & \textbf{(-0.328, -0.172)} & \textbf{(-0.335, -0.215)} \\
 & \textbf{0.14, 0.58} & \textbf{0.212, 0.903} & \textbf{0.156, 0.894} & \textbf{0.121, 0.916} \\
\bottomrule
\multirow{6}{1 cm}{0.25} & 0.005 & 0.004 & 0.001 & 0.001\\
 & (0.022, 0.286) & (0.187, 0.445) & (0.095, 0.245) & (0.173, 0.297)\\
 & 0.265, 0.825 & 0.257, 0.953 & 0.149, 0.873 & 0.124, 0.925 \\
 & \textbf{0.005} & \textbf{0.004} & \textbf{0.001} & \textbf{0.001}\\
 & \textbf{(0.022, 0.286)} & \textbf{(0.197, 0.435)} & \textbf{(0.102, 0.238)} & \textbf{(0.178, 0.292)} \\
 & \textbf{0.263, 0.825} & \textbf{0.237, 0.942} & \textbf{0.137, 0.853} & \textbf{0.115, 0.904} \\
\bottomrule
\multirow{6}{1 cm}{0.5} & 0.011 & 0.009 & 0.002 & 0.002\\
 & (0.117, 0.533) & (0.397, 0.761) & (0.158, 0.332) & (0.4, 0.59)\\
 & 0.416, 0.87 & 0.365, 0.955 & 0.173, 0.819 & 0.19, 0.924 \\
 & \textbf{0.011} & \textbf{0.008} & \textbf{0.002} & \textbf{0.002}\\
 & \textbf{(0.121, 0.529)} & \textbf{(0.409, 0.749)} & \textbf{(0.163, 0.327)} & \textbf{(0.404, 0.586)} \\
 & \textbf{0.409, 0.865} & \textbf{0.341, 0.943} & \textbf{0.163, 0.799} & \textbf{0.183, 0.917} \\
\bottomrule
\multirow{6}{1 cm}{0.75} &  0.023 & 0.012 & 0.007 & 0.004\\
 & (0.351, 0.939) & (0.502, 0.928) & (0.555, 0.887) & (0.586, 0.842)\\
 & 0.588, 0.89 & 0.425, 0.923 & 0.332, 0.933 & 0.256, 0.96 \\
 & \textbf{0.021} & \textbf{0.011} & \textbf{0.006} & \textbf{0.004}\\
 & \textbf{(0.362, 0.928)} & \textbf{(0.51, 0.92)} & \textbf{(0.563, 0.879)} & \textbf{(0.596, 0.832)} \\
 & \textbf{0.565, 0.876} & \textbf{0.411, 0.913} & \textbf{0.316, 0.92} & \textbf{0.236, 0.943} \\
\bottomrule
\end{tabular*}
\label{Table of CI Normality 2}
\end{table}

\begin{table}
\caption{Estimate of the variance of MLE \texorpdfstring{$\hat{\mu}$}{TEXT}, CI of $\mu$ assuming an asymptotic normality of MLE \texorpdfstring{$\hat{\mu}$}{TEXT}, its AL and CP, and all these entities removing outliers for $\mu$, sample size $n$, and simulation size $N=1000$}
\renewcommand{\arraystretch}{1}
\begin{tabular*}{\textwidth}{@{\extracolsep\fill}ccccc} 
\toprule
\multirow{2}{1 cm}{$\mu$} & \multicolumn{4}{@{}c@{}}{n}  \\ \cmidrule{2-5}%
& 50 & 100 & 150 & 250 \\ \midrule
\multirow{6}{1 cm}{1} &  0.021 & 0.019 & 0.01 & 0.007\\
 & (0.303, 0.867) & (0.715, 1.261) & (0.753, 1.137) & (0.827, 1.157)\\
 & 0.564, 0.787 & 0.547, 0.943 & 0.383, 0.887 & 0.33, 0.93 \\
 & \textbf{0.018} & \textbf{0.018} & \textbf{0.009} & \textbf{0.006}\\
 & \textbf{(0.32, 0.85)} & \textbf{(0.723, 1.253)} & \textbf{(0.764, 1.126)} & \textbf{(0.837, 1.147)} \\
 & \textbf{0.531, 0.762} & \textbf{0.53, 0.937} & \textbf{0.362, 0.872} & \textbf{0.31, 0.922} \\
\bottomrule
\multirow{6}{1 cm}{1.5} & 0.088 & 0.026 & 0.021 & 0.014\\
 & (0.797, 1.961) & (0.867, 1.497) & (1.068, 1.638) & (1.236, 1.706)\\
 & 1.164, 0.931 & 0.631, 0.845 & 0.569, 0.891 & 0.471, 0.922 \\
 & \textbf{0.082} & \textbf{0.024} & \textbf{0.019} & \textbf{0.013}\\
 & \textbf{(0.817, 1.941)} & \textbf{(0.881, 1.483)} & \textbf{(1.082, 1.624)} & \textbf{(1.247, 1.695)} \\
 & \textbf{1.124, 0.92} & \textbf{0.602, 0.821} & \textbf{0.541,0.878} & \textbf{0.448, 0.909} \\
\bottomrule
\multirow{6}{1 cm}{2} &  0.33 & 0.125 & 0.034 & 0.045\\
 & (1.159, 3.411) & (1.544, 2.93) & (1.309, 2.031) & (1.876, 2.71)\\
 & 2.253, 0.983 & 1.386, 0.98 & 0.722, 0.877 & 0.834, 0.976 \\
 & \textbf{0.279} & \textbf{0.109} & \textbf{0.032} & \textbf{0.041}\\
 & \textbf{(1.251, 3.319)} & \textbf{(1.59, 2.884)} & \textbf{(1.318, 2.022)} & \textbf{(1.898, 2.688)} \\
 & \textbf{2.069, 0.973} & \textbf{1.294, 0.971} & \textbf{0.704, 0.866} & \textbf{0.79, 0.966} \\
\bottomrule
\multirow{6}{1 cm}{3} & 9.448 & 0.604 & 0.11 & 0.06\\
 & (-2.802, 9.246) & (1.9, 4.946) & (1.876, 3.178) & (1.981, 2.941)\\
 & 12.049, 0.999 & 3.046, 0.982 & 1.303, 0.846 & 0.961, 0.829 \\
 & \textbf{0.747} & \textbf{0.453} & \textbf{0.098} & \textbf{0.056}\\
 & \textbf{(1.528, 4.916)} & \textbf{(2.104, 4.742)} & \textbf{(1.912, 3.142)} & \textbf{(1.998, 2.924)} \\
 & \textbf{3.387, 0.944} & \textbf{2.638, 0.971} & \textbf{1.23, 0.815} & \textbf{0.926, 0.812} \\
\bottomrule
\multirow{6}{1 cm}{4} &  144.008 & 1944.417 & 1.765 & 0.934\\
 & (-19.251, 27.789) & (-79.605, 93.247) & (1.931, 7.139) & (2.97, 6.758)\\
 & 47.04, 0.992 & 172.851, 1 & 5.207, 0.988 & 3.788, 0.989 \\
 & \textbf{1.803} & \textbf{5.673} & \textbf{0.808} & \textbf{0.61}\\
 & \textbf{(1.637, 6.901)} & \textbf{(2.153, 11.489)} & \textbf{(2.773, 6.297)} & \textbf{(3.334, 6.394)} \\
 & \textbf{5.263, 0.904} & \textbf{9.337, 0.992} & \textbf{3.524, 0.963} & \textbf{3.061, 0.982} \\
\bottomrule
\multirow{6}{1 cm}{5} & 31.615 & 568.554 & 0.368 & 0.751 \\
 & (-7.124, 14.916) & (-41.835, 51.633) & (2.278, 4.656) & (2.881, 6.277)\\
 & 22.041, 0.951 & 93.468, 0.996 & 2.378, 0.667 & 3.396, 0.92 \\
 & \textbf{1.392} & \textbf{1.58} & \textbf{0.27} & \textbf{0.45}\\
 & \textbf{(1.584, 6.208)} & \textbf{(2.435, 7.363)} & \textbf{(2.448, 4.486)} & \textbf{(3.264, 5.894)} \\
 & \textbf{4.624, 0.722} & \textbf{4.928, 0.885} & \textbf{2.038, 0.594} & \textbf{2.631, 0.843} \\
\bottomrule
\end{tabular*}
\label{Table of CI Normality 3}
\end{table}

\subsection{Using the empirical distribution of the MLE \texorpdfstring{$\hat{\mu}$}{TEXT}}
\label{Sec. construction of CI using empirical distribution of MLE}

Here, we construct the CI for $\mu$ using the empirical distribution of the MLE $\hat{\mu}$ based on simulation instead of assuming asymptotic normality for $\hat{\mu}$. Following Algorithm \ref{algo:CI using empirical distr of MLE} and Algorithm \ref{algo:AL and CP for CI using empirical distr of MLE} gives the procedure of construction of CI, and its AL and CP respectively. 

\begin{algorithm}
    \caption{To compute CI for $\mu$ using empirical distribution of MLE $\hat{\mu}$}
    \label{algo:CI using empirical distr of MLE}
    Perform Steps 1-3 of Algorithm \ref{algo:CI using asymp normality of MLE}.\\
    Arrange $N$ values of $\hat{\mu}$ obtained in Step 1 in ascending order.\\
    The lower ($L$) and upper ($U$) limit of $100(1-\alpha)\%$ CI for $\mu$ is 
    $$
    \begin{array}{ll}
         L & = \left( \frac{N\alpha}{2} \right)^{th} \text{ordered value of } \hat{\mu} \\
        U & =  \left[ N\left (1 - \frac{\alpha}{2} \right) \right]^{th} \text{ordered value of } \hat{\mu}.
    \end{array}
    $$   
\end{algorithm}

\begin{algorithm}
    \caption{To compute the AL and CP of the CI for $\mu$ constructed using the empirical distribution of MLE $\hat{\mu}$}
    \label{algo:AL and CP for CI using empirical distr of MLE}
    Generate a random sample of size $n$ from $SSLUD(\mu)$ for a specified value of $\mu$.\\
    Perform all Steps of Algorithm \ref{algo:CI using empirical distr of MLE} and compute length of CI = $U-L$.\\
    Repeat above Steps 1-2 for $N$ times.\\
    The AL of these CI is the average length of $N$ CIs obtained in Step 3 and the corresponding CP is the proportion of CIs obtained in Step 3 which contain $\mu$ specified in Step 1.
\end{algorithm}

The modified CI based on this method is constructed by removing the outliers using the IQR method of outlier detection as discussed in Section \ref{Sec. construction of CI using asymptotic normal distribution of the MLE}. To compute this modified CI one has to perform Steps 1-2 of Algorithm \ref{algo:CI using empirical distr of MLE} and remove outliers out of $N$ values of $\hat{\mu}$ using the IQR method. Suppose $N^*$ denotes the number of remaining values of $\hat{\mu}$ after removing outliers ($N^* \leqslant N$). Now obtain modified CI for $\mu$ by replacing $N$ with $N^*$ using Step 3 of Algorithm \ref{algo:CI using empirical distr of MLE} so that ordered $N^*$ number of $\hat{\mu}$ values are used. The AL and CP of these modified CI are computed using Algorithm \ref{algo:AL and CP for CI using empirical distr of MLE} by applying the modification discussed above in Step 2.

Table \ref{Tab CI using empirical 1} and Table \ref{Tab CI using empirical 2} represent the CI of $\mu$ using empirical distribution of MLE, its AL and CP, and all these entities removing outliers for $\mu=-5, -4, -3, -2, -1.5, -1, -0.75, -0.5, -0.25, 0.25, 0.5, 0.75, 1, 1.5, 2, 3, 4, 5$ and the sample size $n= 50, 100, 150, 250$ using simulation size $N=1000$. From these tables it is seen that the CI constructed using empirical distribution of the MLE $\hat{\mu}$ is overall comparable with the CI constructed assuming asymptotic normality of the MLE $\hat{\mu}$. However, in some cases, this method performs worst in terms of CP (CP=0). The same is observed for its modified CI also. Hence, one can say that the method described in Section \ref{Sec. construction of CI using asymptotic normal distribution of the MLE} is more reliable than the method described in Section \ref{Sec. construction of CI using empirical distribution of MLE}. 

\begin{table}
\caption{
CI of $\mu$ using the empirical distribution of MLE \texorpdfstring{$\hat{\mu}$}{TEXT}, its AL and CP, and all these entities removing outliers for $\mu$, sample size $n$, and simulation size $N=1000$} \label{Tab CI using empirical 1}

\renewcommand{\arraystretch}{1}
\begin{tabular*}{\textwidth}{@{\extracolsep\fill}ccccc} 
\toprule
\multirow{2}{1 cm}{$\mu$} & \multicolumn{4}{@{}c@{}}{n}  \\ \cmidrule{2-5}%
& 50 & 100 & 150 & 250 \\ \midrule
\multirow{4}{1 cm}{-5} & (-64.466, 32.803) & (-13.737, -2.869) & (-7.735, -2.951) & (-11.658, -4.048)\\
 & 93.947, 1 & 10.843, 1 & 4.889, 1 & 6.779, 1 \\
 & \textbf{(-14.912, -2.521)} & \textbf{(-8.749, -2.869)} & \textbf{(-6.647, -2.94)} & \textbf{(-8.788, -4.04)} \\
 & \textbf{12.462, 1} & \textbf{5.715, 1} & \textbf{3.631, 1} & \textbf{4.84, 1} \\
\bottomrule
\multirow{4}{1 cm}{-4} & (-2.484, -1.01) & (-4.589, -2.032) & (-13.31, -3.354) & (-5.155, -2.913)\\
 & 1.523, 0 & 2.634, 1  & 11.056, 1 & 2.232, 1 \\
 & \textbf{(-2.461, -1.013)} & \textbf{(-4.359, -2.032)} & \textbf{(-9.901, -3.352)} & \textbf{(-4.93, -2.917)} \\
 & \textbf{1.429, 0} & \textbf{2.273, 1} & \textbf{6.209, 1} & \textbf{2, 1} \\
\bottomrule
\multirow{4}{1 cm}{-3} & (-6.083, -1.674) & (-3.955, -1.8) & (-4.005, -2.131) & (-3.586, -2.247)\\
 & 4.458, 1 & 1.91, 1 & 1.853, 1 & 1.289, 1 \\
 & \textbf{(-4.849, -1.675)} & \textbf{(-3.578, -1.801)} & \textbf{(-3.916, -2.146)} & \textbf{(-3.537, -2.253)} \\
 & \textbf{3.214, 1} & \textbf{1.744, 1} & \textbf{1.695, 1} & \textbf{1.216, 1} \\
\bottomrule
\multirow{4}{1 cm}{-2} & (-2.701, -1.022) & (-1.84, -1.006) & (-2.366, -1.426) & (-1.966, -1.384)\\
 & 1.652, 1 & 0.809, 0 & 0.924, 1 & 0.575, 0.001 \\
 & \textbf{(-2.575, -1.032)} & \textbf{(-1.82, -1.007)} & \textbf{(-2.341, -1.461)} & \textbf{(-1.943, -1.408)} \\
 & \textbf{1.539, 1} & \textbf{0.775, 0} & \textbf{0.882, 1} & \textbf{0.549, 0} \\
\bottomrule
\multirow{4}{1 cm}{-1.5} & (-2.071, -0.845) & (-1.826, -1.008) & (-1.704, -1.079) & (-1.56, -1.13)\\
 & 1.248, 1 & 0.797, 1 & 0.617, 1 & 0.448, 1 \\
 & \textbf{(-2.049, -0.855)} & \textbf{(-1.792, -1.025)} & \textbf{(-1.693, -1.1)} & \textbf{(-1.559, -1.152)} \\
 & \textbf{1.186, 1} & \textbf{0.765, 1} & \textbf{0.59, 1} & \textbf{0.428, 1} \\
\bottomrule
\multirow{4}{1 cm}{-1} & (-1.571, -0.565) & (-1.163, -0.619) & (-0.967, -0.591) & (-1.013, -0.713)\\
 & 0.949, 1 & 0.531, 1 & 0.37, 0 & 0.296, 0.999 \\
 & \textbf{(-1.543, -0.58)} & \textbf{(-1.159, -0.658)} & \textbf{(-0.967, -0.6)} & \textbf{(-1.013, -0.719)} \\
 & \textbf{0.908, 1} & \textbf{0.508, 1} & \textbf{0.353, 0} & \textbf{0.282, 0.999} \\
\bottomrule
\multirow{4}{1 cm}{-0.75} & (-1.01, -0.372) & (-0.79, -0.399) & (-0.87, -0.533) & (-0.812, -0.552)\\
 & 0.68, 1 & 0.39, 1 & 0.338, 1 & 0.245, 1 \\
 & \textbf{(-1.001, -0.391)} & \textbf{(-0.786, -0.407)} & \textbf{(-0.869, -0.541)} & \textbf{(-0.807, -0.565)} \\
 & \textbf{0.648, 1} & \textbf{0.371, 1} & \textbf{0.322, 1} & \textbf{0.232, 1} \\
\bottomrule
\multirow{4}{1 cm}{-0.5} & (-0.424, -0.037) & (-0.68, -0.325) & (-0.498, -0.266) & (-0.497, -0.326)\\
 & 0.388, 0 & 0.349, 1 & 0.231, 0.264 & 0.176, 0.872 \\
 & \textbf{(-0.422, -0.037)} & \textbf{(-0.675, -0.347)} & \textbf{(-0.498, -0.274)} & \textbf{(-0.497, -0.337)} \\
 & \textbf{0.381, 0} & \textbf{0.331, 1} & \textbf{0.217, 0.172} & \textbf{0.166, 0.816} \\
\bottomrule
\multirow{4}{1 cm}{-0.25} & (-0.131, -0.001) & (-0.342, -0.112) & (-0.306, -0.142) & (-0.322, -0.189)\\
 & 0.13, 0 & 0.234, 1 & 0.176, 1 & 0.135, 1 \\
 & \textbf{(-0.128, -0.001)} & \textbf{(-0.341, -0.127)} & \textbf{(-0.306, -0.149)} & \textbf{(-0.321, -0.198)} \\
 & \textbf{0.127, 0} & \textbf{0.217, 1} & \textbf{0.164, 1} & \textbf{0.127, 1} \\
\bottomrule
\end{tabular*}
\label{Table of }
\end{table}

\begin{table}
\caption{
CI of $\mu$ using the empirical distribution of MLE \texorpdfstring{$\hat{\mu}$}{TEXT}, its AL and CP, and all these entities removing outliers for $\mu$, sample size $n$, and simulation size $N=1000$}
\label{Tab CI using empirical 2}
\renewcommand{\arraystretch}{1}
\begin{tabular*}{\textwidth}{@{\extracolsep\fill}ccccc} 
\toprule
\multirow{2}{1 cm}{$\mu$} & \multicolumn{4}{@{}c@{}}{n}  \\ \cmidrule{2-5}%
& 50 & 100 & 150 & 250 \\ \midrule
\multirow{4}{1 cm}{0.25} & (0.003, 0.244) & (0.148, 0.401) & (0.069, 0.218) & (0.157, 0.277)\\
 & 0.236, 0.022 & 0.256, 1 & 0.148, 0 & 0.124, 1 \\
 & \textbf{(0.003, 0.241)} & \textbf{(0.166, 0.398)} & \textbf{(0.078, 0.218)} & \textbf{(0.165, 0.276)} \\
 & \textbf{0.235, 0.013} & \textbf{0.238, 1} & \textbf{0.137, 0} & \textbf{0.117, 1} \\
\bottomrule
\multirow{4}{1 cm}{0.5} & (0.036, 0.469) & (0.355, 0.727) & (0.132, 0.303) & (0.376, 0.57)\\
 & 0.412, 0 & 0.362, 1 & 0.174, 0 & 0.19, 1 \\
 & \textbf{(0.035, 0.463)} & \textbf{(0.375, 0.72)} & \textbf{(0.139, 0.304)} & \textbf{(0.381, 0.57)} \\
 & \textbf{0.389, 0} & \textbf{0.343, 1} & \textbf{0.162, 0} & \textbf{0.179, 1} \\
\bottomrule
\multirow{4}{1 cm}{0.75} & (0.243, 0.846) & (0.445, 0.883) & (0.5, 0.848) & (0.553, 0.814)\\
 & 0.593, 1 & 0.419, 1 & 0.332, 1 & 0.245, 1 \\
 & \textbf{(0.272, 0.837)} & \textbf{(0.455, 0.88)} & \textbf{(0.515, 0.845)} & \textbf{(0.565, 0.806)} \\
 & \textbf{0.563, 1} & \textbf{0.399, 1} & \textbf{0.316, 1} & \textbf{0.232, 1} \\
\bottomrule
\multirow{4}{1 cm}{1} & (0.214, 0.782) & (0.671, 1.187) & (0.711, 1.105) & (0.796, 1.107)\\
 & 0.556, 0 & 0.543, 1 & 0.413, 1 & 0.321, 1 \\
 & \textbf{(0.247, 0.782)} & \textbf{(0.682, 1.187)} & \textbf{(0.727, 1.09)} & \textbf{(0.808, 1.107)} \\
 & \textbf{0.527, 0} & \textbf{0.519, 1} & \textbf{0.393, 1} & \textbf{0.306, 1} \\
\bottomrule
\multirow{4}{1 cm}{1.5} & (0.743, 1.87) & (0.787, 1.427) & (1.023, 1.591) & (1.203, 1.665)\\
 & 1.123, 1 & 0.641, 0 & 0.583, 1 & 0.478, 1 \\
 & \textbf{(0.747, 1.847)} & \textbf{(0.793, 1.416)} & \textbf{(1.04, 1.585)} & \textbf{(1.219, 1.664)} \\
 & \textbf{1.072, 1} & \textbf{0.614, 0} & \textbf{0.559, 1} & \textbf{0.457, 1} \\
\bottomrule
\multirow{4}{1 cm}{2} & (1.234, 3.391) & (1.513, 2.915) & (1.254, 1.992) & (1.853, 2.698)\\
 & 2.192, 1 & 1.387, 1 & 0.74, 0.191 & 0.85, 1 \\
 & \textbf{(1.234, 3.22)} & \textbf{(1.518, 2.835)} & \textbf{(1.259, 1.979)} & \textbf{(1.874, 2.667)} \\
 & \textbf{1.964, 1} & \textbf{1.303, 1} & \textbf{0.709, 0.023} & \textbf{0.812, 1} \\
\bottomrule
\multirow{4}{1 cm}{3} & (1.617, 6.216) & (2.127, 5.177) & (1.828, 3.124) & (1.93, 2.904)\\
 & 4.617, 1 & 3.051, 1 & 1.322, 1 & 0.948, 0 \\
 & \textbf{(1.615, 5.019)} & \textbf{(2.117, 4.748)} & \textbf{(1.849, 3.091)} & \textbf{(1.938, 2.858)} \\
 & \textbf{3.278, 1} & \textbf{2.553, 1} & \textbf{1.244, 1} & \textbf{0.903, 0} \\
\bottomrule
\multirow{4}{1 cm}{4} & (1.972, 13.352) & (3.205, 33.755) & (3.008, 7.302) & (3.491, 7.136)\\
 & 11.734, 1 & 29.377, 1 & 4.662, 1 & 3.848, 1 \\
 & \textbf{(1.991, 7.117)} & \textbf{(3.483, 12.964)} & \textbf{(3.006, 6.548)} & \textbf{(3.49, 6.555)} \\
 & \textbf{5.35, 1} & \textbf{9.984, 1} & \textbf{3.521, 1} & \textbf{3.139, 1} \\
\bottomrule
\multirow{4}{1 cm}{5} & (1.857, 9.654) & (2.772, 10.66) & (2.425, 4.72) & (3.274, 6.5)\\
 & 8.474, 1 & 8.718, 1 & 2.388, 0.019 & 3.288, 1 \\
 & \textbf{(1.867, 6.519)} & \textbf{(2.774, 7.73)} & \textbf{(2.424, 4.427)} & \textbf{(3.27, 5.941)} \\
 & \textbf{4.522, 1} & \textbf{5.043, 1} & \textbf{2.109, 0} & \textbf{2.766, 1} \\
\bottomrule
\end{tabular*}
\label{Table of Bias and MSE 1}
\end{table}

\section{Application} \label{sec application}
In this section, we illustrate the construction of CI for $\mu$ of $SSLUD(\mu)$ discussed in Section \ref{Sec. Construction of confidence interval} using data of transformed daily percentage change $Z_t$ in the price of NIFTY 50, an Indian stock market index of size $n=82$. It is referred in \cite{Lohot2023Dixit} and validated that $SSLUD(\mu)$ is suitable for modeling $Z_t$. The transformed data values of $Z_t$ are as follows.

 \noindent -\,0.64, -\,2.33, -\,2.98, 0.14, 0.30, -\,0.11, -\,1.20, -\,0.31, 0.06, -\,0.91, -\,0.86, 0.07, 0.77, 0.22, -\,0.13, -\,1.80, -\,0.42, 0.27, -\,0.51, 0.07, -\,0.55, -\,0.81, -\,0.51, -\,1.87, -\,1.76, -\,1.81, -\,1.59, -\,3.46, -\,0.05, -\,1.77, -\,0.85, 0.59, 0.57, 0.36, -\,2.04, -\,1.05, -\,2.53, -\,0.49, 0.34, 0.01, -\,2.11, -\,3.86, 2.23, -\,0.97, -\,0.90, -\,0.96, -\,1.20, -\,1.47, -\,0.97, -\,5.58, 1.73, 0.01, -\,1.92, -\,1.45, -\,2.33, -\,3.15, 0.15, 1.27,  0.73, -\,0.59, 0.65, -\,2.03, 1.07, 1.04, -\,1.78, 0.36, -\,1.20, -\,0.93, -\,1.20, -\,0.40, -\,0.20, 0.20, -\,0.99, 0.38, 1.37, -\,1.33, -\,1.63, -\,1.74, 0.02, -\,1.42, -\,1.62, -\,1.11.

The MLE $\hat{\mu}$ of $\mu$ for this data is -2.589. The estimate of the variance of $\hat{\mu}$ and the corresponding CI for $\mu$ assuming the asymptotic normality of $\hat{\mu}$ as discussed in Section \ref{Sec. construction of CI using asymptotic normal distribution of the MLE} are 0.2698767 and (-3.607195, -1.570805) respectively. Further, its modified values are 0.1927435 and (-3.449474, -1.728526) respectively after removing the outliers. The CI for $\mu$ using the empirical distribution of the MLE $\hat{\mu}$ as discussed in Section \ref{Sec. construction of CI using empirical distribution of MLE} is (-3.760858, -1.602701) and its modified CI is (-3.420354, -1.622304) after removing the outliers. Overall both methods have given reasonably the same CIs.

\section*{Declarations}

\subsection*{Conflict of interest} The authors declare that they have no conflict of interest to disclose.

\subsection*{Ethics approval and consent to participate} Not applicable.

\subsection*{Consent for publication} Both authors have agreed to submit and publish this paper.

\clearpage
\bibliography{citation}

\begin{thebibliography}{}

\bibitem[Azzalini, 1985]{azzalini1985class}
Azzalini, A. (1985).
\newblock A class of distributions which includes the normal ones.
\newblock {\em Scandinavian journal of statistics}, 12(2):171--178.

\bibitem[{Lohot, R. K. and Dixit, V. U.}, 2024]{Lohot2023Dixit}
{Lohot, R. K. and Dixit, V. U.} (2024).
\newblock {The skew-symmetric-Laplace-uniform distribution}.
\newblock \url{https://arxiv.org/abs/2406.10805}.
\newblock arXiv:2406.10805 [math.ST].

\end{thebibliography}
\bibliographystyle{apalike}

\end{document}